\documentclass[reqno]{article}
\usepackage{tikz}
\usepackage{amsmath,amsthm,amssymb}
\usepackage{graphicx}
\graphicspath{ {./Steiner/} }
\newtheorem{theorem}{Theorem}[section]
\newtheorem{lemma}[theorem]{Lemma}
\newtheorem{proposition}[theorem]{Proposition}
\newtheorem{corollary}[theorem]{Corollary}
\newtheorem{conjecture}[theorem]{Conjecture}

\newtheorem{definition}[theorem]{Definition}

\def\beq{\begin{equation}}\def\eeq{\end{equation}}
\def\beqn{\begin{eqnarray}}\def\eeqn{\end{eqnarray}}

\def\eps{\varepsilon}
\def\qed{\ifhmode\unskip\nobreak\fi\quad\ifmmode\Box\else$\Box$\fi}

\begin{document}
\title{Tur\'an and Ramsey numbers in linear triple systems}
\author{Andr\'as Gy\'arf\'as\thanks{Alfr\'ed R\'enyi Institute of Mathematics, Budapest, P.O. Box 127, Budapest, Hungary, H-1364. \texttt{gyarfas.andras@renyi.hu}, \texttt{sarkozy.gabor@renyi.hu}} \thanks{Research supported in part by
NKFIH Grant No. K116769.} \and G\'{a}bor N. S\'ark\"ozy\footnotemark[1]
\thanks{Computer Science Department, Worcester Polytechnic Institute, Worcester, MA.} \thanks{Research supported in part by
NKFIH Grants No. K116769, K117879.}
}

\maketitle

\begin{abstract}
In this paper we study Tur\'an and Ramsey numbers in linear triple systems, defined as $3$-uniform
hypergraphs in which any two triples intersect in at most one
vertex.

A famous result of Ruzsa and Szemer\'edi  is that for any fixed $c>0$ and large enough $n$ the following Tur\'an-type theorem holds.
If a linear triple system on $n$ vertices has at least
$cn^2$ edges then it contains a {\em triangle}: three pairwise intersecting triples without a common vertex. In
this paper we extend this result from triangles to other triple systems, called {\em
$s$-configurations}. The main tool is a generalization of the induced matching lemma from $aba$-patterns to more general ones.

We slightly generalize $s$-configurations to {\em extended
$s$-configurations}. For these we cannot prove the corresponding Tur\'an-type theorem, but we prove that they have the weaker, Ramsey property: they
can be found in any $t$-coloring of the blocks of any sufficiently
large Steiner triple system. Using this, we show that all
unavoidable configurations with at most 5 blocks, except possibly
the ones containing the sail $C_{15}$ (configuration with blocks
123, 345, 561 and 147), are $t$-Ramsey for any $t\geq 1$. The most interesting one among them is the {\em wicket}, $D_4$,
formed by three rows and two columns of a $3\times 3$ point matrix.
In fact, the wicket is $1$-Ramsey in a very strong sense: all Steiner triple systems except the Fano plane must contain a wicket.

\end{abstract}

\section{Introduction}

\subsection{Tur\'an-type problems}\label{subturan}
In this paper we study {\em linear triple systems}, defined as $3$-uniform
hypergraphs in which any two triples intersect in at most one
vertex. The $(k,\ell)$-family is the family of all linear triple systems with $\ell$ triples on at most $k$ vertices. A
famous conjecture of Brown, Erd\H{o}s and T. S\'{o}s \cite{BETS}
claims the following. (It is well-known that here we can restrict our attention to linear triple systems, and that is why
in this paper we focus only on linear triple systems.)
\begin{conjecture}\label{bes} If a linear triple system on $n$ vertices does not contain any member of the $(k+3,k)$-family then it has $o(n^2)$ triples.
\end{conjecture}

For $k=3$ there is only one member in the $(6,3)$-family, the {\em triangle}, three pairwise intersecting triples without a common vertex.
Conjecture \ref{bes} in this case was famously proved by Ruzsa and
Szemer\'edi \cite{RS}, in addition with the surprising lower bound: there are triangle-free linear triple systems with $n^{2-o(1)}$ triples. This became known as the (6,3)-theorem.
The (6,3)-theorem had a huge influence. For example the celebrated Triangle Removal
Lemma (see \cite{CF} for a survey) was devised in order to find another proof
for the (6,3)-theorem.

Another simple proof of the (6,3)-theorem was found later
by Szemer\'edi (see \cite{KS}). In this proof Szemer\'edi used the
argument that if a dense graph is properly edge colored with $O(n)$ colors then
there is a path on 3 edges that gets only 2 colors (called an $aba$-pattern).
This (or its contrapositive) became known as the ``induced matching lemma".
S\'ark\"ozy and Selkow \cite{SS}  generalized this argument; under the same conditions
certain subtrees that get few colors can be found.
This approach led to the following result, showing that Conjecture \ref{bes} is close to being true.

\begin{theorem}[\cite{SS}]\label{stan} If a linear triple system on $n$ vertices does not contain any member of
the $(k+2 +\lfloor \log_2 k \rfloor,k)$-family then it has  $o(n^2)$
triples.
\end{theorem}

However, it still remained open whether one can replace the $\lfloor
\log_2 k \rfloor$ term with 1 and prove Conjecture \ref{bes}. For
example the $k=4$ case of Conjecture \ref{bes}, the $(7,4)$-problem
is still open. The $(7,4)$-family has three members (see Figure 4) and the problem boils down to deciding whether one of $C_{14}$
and $C_{16}$ must be present in any linear triple system on $n$ vertices with at least $cn^2$ edges, see \cite{FGY1}.

What we do get from Theorem \ref{stan} is that
forbidding $(8,4),(9,5)$, etc.  families, we have $o(n^2)$ triples.
Since then the only progress was obtained by Solymosi and Solymosi
\cite{SOSO} who improved the $k=10$ case of Theorem \ref{stan}
showing that forbidding $(14,10)$-families already implies $o(n^2)$
triples (instead of the $(15,10)$ implied by Theorem \ref{stan}).
Recently, Conlon, Gishboliner, Levanzov and Shapira \cite{CGLS}
announced an improvement of Theorem \ref{stan} in which the
$\log{k}$ term was replaced by a $O(\log{k}/\log{\log{k}})$.

In this paper we refine and generalize further the argument that led
to Theorem \ref{stan}. In fact, we will define a special class of
linear triple systems, called $s$-configurations (not to be confused with
$(k,\ell)$-configurations), and we will show that these must occur
in every linear triple system on $n$ vertices with at least $cn^2$ triples,
i.e. can play the role of the $(k+3,k)$-family in Conjecture
\ref{bes} (or the role
of the family in Theorem \ref{stan}). The main feature here is that we look for {\em one specific $(k,\ell)$-configuration}
 instead of {\em a family}.
For example, we will show that any of the $(9,5)$-configurations $D_1,D_2$ and $D_3$  in Figures 1 and 2 must be present in linear triple systems on $n$
vertices with at least $cn^2$ triples. Note that for $D_4$ in Figure 2, called the {\em wicket}, we are able to prove only a weaker,
Ramsey-type result (see in the next subsection). However, extending
the wicket to the {\em grid} defined by three rows and three columns
of a $3\times 3$ point matrix, F\"uredi and Ruszink\'o conjectured
\cite{FR} that grid-free linear triple systems on $n$ points can have $cn^2$ (in fact even ${n(n-1)\over 6}$) triples.
The weaker conjecture was proved recently with $c={1\over 16}$ in \cite{GS}.

\begin{figure}[ht]
\centering
\includegraphics[scale=.7]{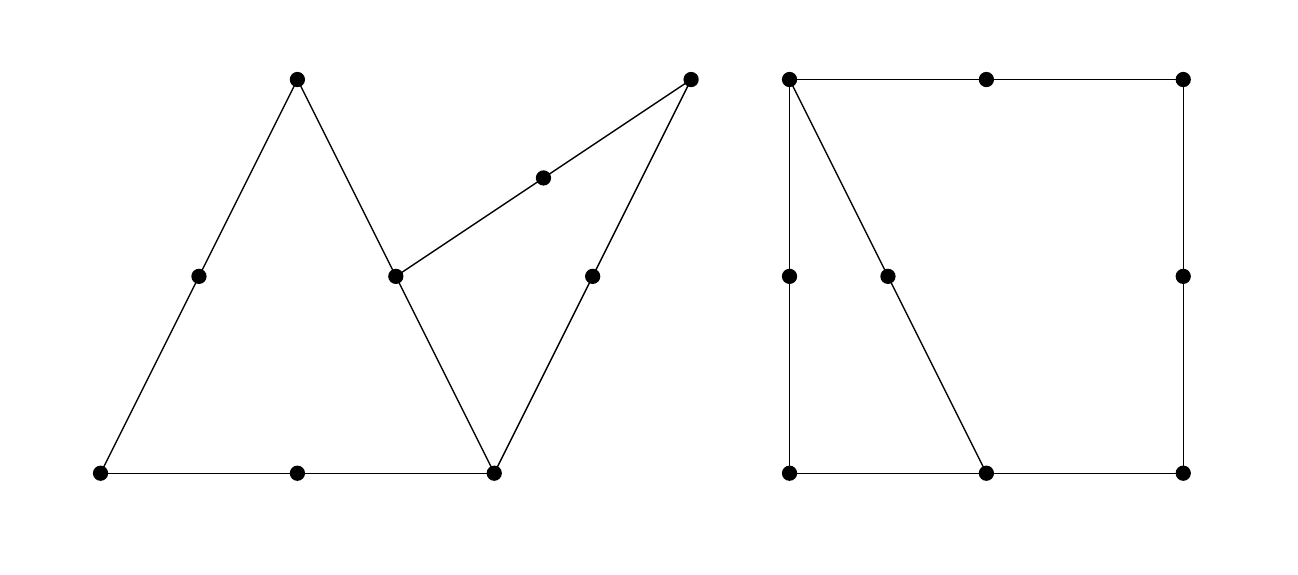}
\caption{Configurations $D_1$ and $D_2$}
\end{figure}

\begin{figure}[ht]
\centering
\includegraphics[scale=.7]{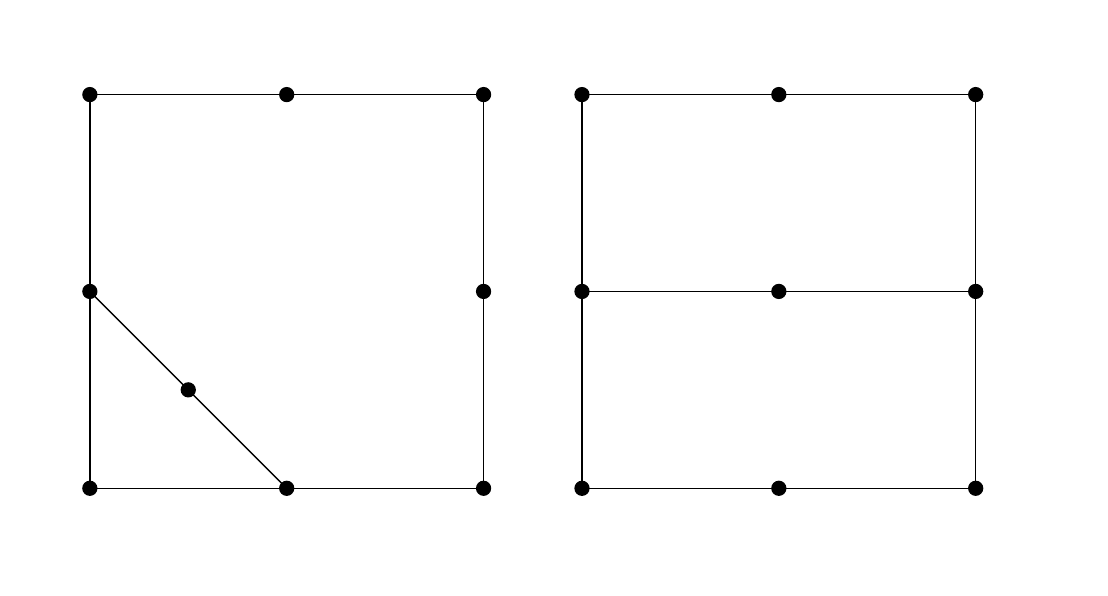}
\caption{Configurations $D_3$ and $D_4$ (wicket)}
\end{figure}

Consider those properly edge-colored forests (acyclic graphs) that can be obtained as the union of  $s$ monochromatic matchings $M^1, \ldots , M^s$
with the following property:
for any $1\le i\le s$ every edge in $M_i$ has a vertex that is not covered by any edge of any $M_j, i < j \leq s$.
We call a forest obtained this way an {\em $s$-pattern}. For example, the path $aba$ is a $2$-pattern but the path $abab$ is not. Note that every
forest has a proper coloring (for some $s$) that makes it an $s$-pattern. Indeed, the coloring where all edges are colored with different colors
is an $s$-pattern for a suitable order of the edges. Next we go one step further: we call a properly edge-colored forest an {\em extended $s$-pattern}
if it is obtained from a disconnected $s$-pattern by joining two of its connected components with a single  edge matching $M^*$ of a new color.
Note that an extended $s$-pattern may or may not be a $(s+1)$-pattern. (The $(s+1)$-pattern is preferable because it leads to stronger results.)
For example, the path $abcab$ is an extended 2-pattern but not a 3-pattern (see Figure 3).

\begin{figure}[ht]
\centering
\includegraphics[scale=.8]{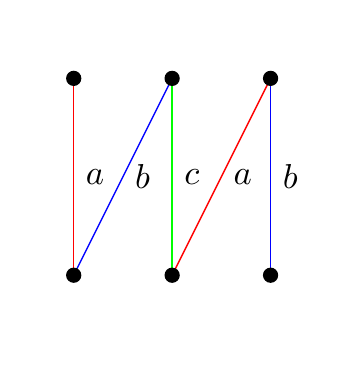}
\caption{The $abcab$ extended 2-pattern}
\end{figure}

A linear triple system ${\cal H}$ is called an {\em
$s$-configuration} (or an {\em extended $s$-configuration}) if it
comes from an $s$-pattern (extended $s$-pattern) by augmenting all
edges $e\in M^i$ with a new {\em augmenting point} $v_i$  to a
triple $e\cup v_i$ in such a way that the $v_i$'s are all distinct
and disjoint from the vertices of the $s$-pattern as well (for
extended $s$-patterns $M^*$ is also augmented with a point that is
distinct from all other augmenting points). For example Figure 4 shows how the
wicket can be obtained by augmenting the extended $2$-pattern $abcab$. (Augmenting
points are shown on the figures by capitalizing the letters of the
corresponding patterns.)

\begin{figure}[ht]
\centering
\includegraphics[scale=.8]{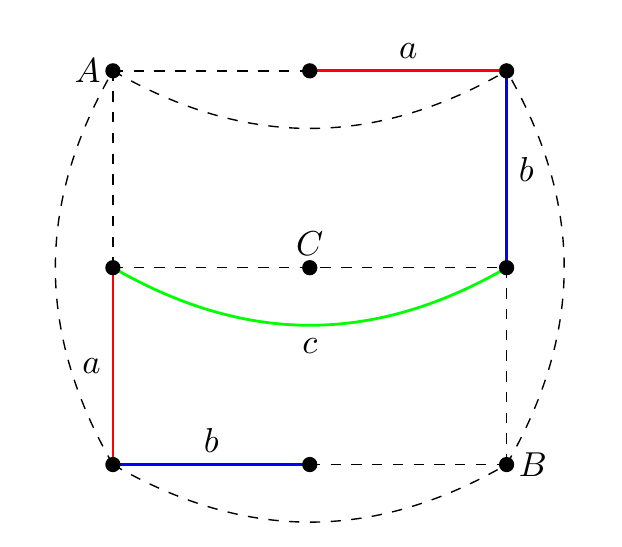}
\caption{Augmenting the $abcab$ extended 2-pattern to a wicket}
\end{figure}

It is worth mentioning that different
patterns may correspond to the same configuration.
For example, $D_3$ is a $3$-configuration from the 3-pattern $abcba$ and also an extended $2$-configuration from the
extended 2-pattern $abacb$ (see Figure 5).

\begin{figure}[ht]
\centering
\includegraphics[scale=.8]{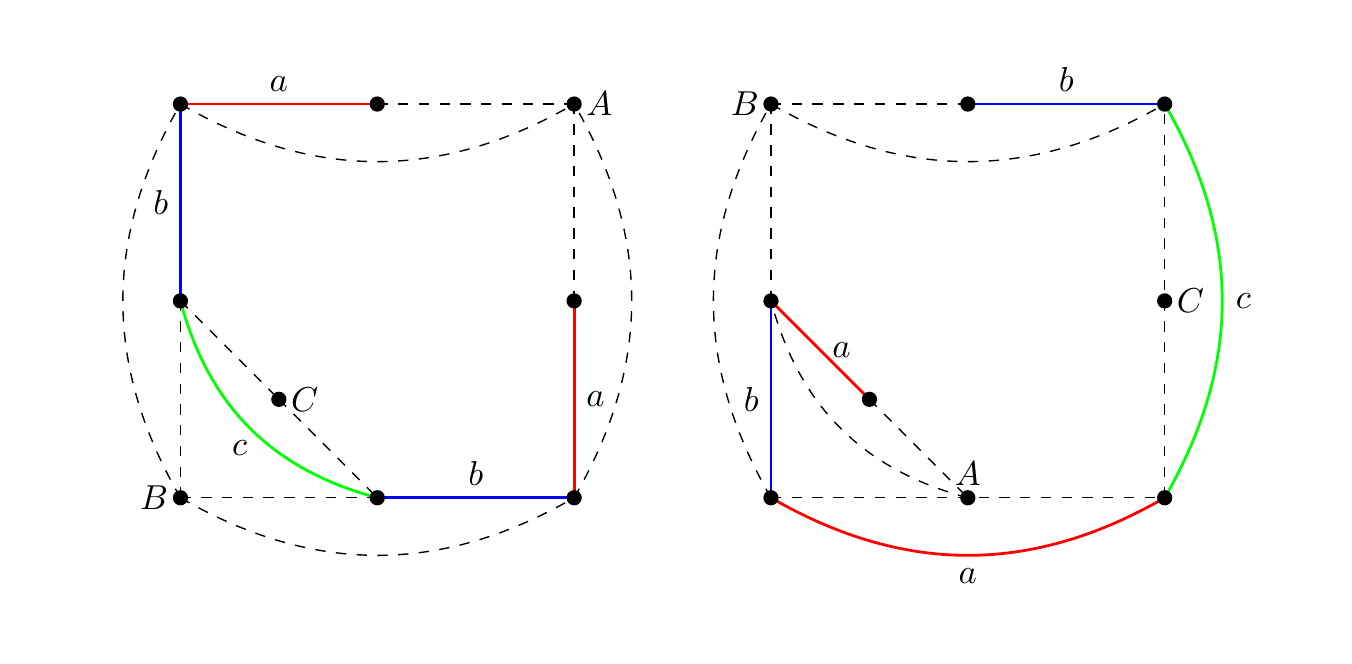}
\caption{Getting $D_3$ from two different patterns}
\end{figure}



\begin{theorem}\label{tetel1}
For every $\delta >0$, integer $s\geq 1$ and $s$-configuration ${\cal H}$,
there is an $n_0=n_0(\delta , s, {\cal H})$ such that if $n\geq n_0$, then every triple system $\cal{G}$ on $n$ vertices with at least $\delta n^2$
triples contains ${\cal H}$ as a subsystem.
\end{theorem}

For example, it is not hard to see that the configurations $D_1,
D_2$ and $D_3$ in Figures 1 and 2 are all 2- or 3-configurations
(see the details in Section 4), and thus these can indeed be found
in linear triple systems on $n$ vertices with at least $\delta n^2$
triples, as claimed above. Note that we do not know this for $D_4$
(since it is an extended 2-configuration).

Theorem \ref{tetel1} will follow from the following generalization of the induced matching lemma.
We say that a properly colored graph $G$ contains a copy of a properly colored graph $H$ if $H$ has a color-preserving isomorphism to a subgraph of $G$.

\begin{theorem}\label{spattern}
For every $\delta,c >0$, integer $s\geq 1$ and $s$-pattern $H$,
there is an $n_0=n_0(\delta , c, s, H)$ with the following property. Let $G$ be a
graph on $n\geq n_0$ vertices with at least $\delta n^2$ edges
that is properly colored by at most $cn$ colors. Then $G$
contains a copy of $H$.
\end{theorem}

Observe that for the $2$-pattern $aba$, Theorem \ref{spattern}  is the induced matching lemma.
But Theorem \ref{spattern} contains many other patterns, for example the pattern used by Duke and R\"odl (Lemma 2 in \cite{DR}). Another example is the pattern $abcba$ which corresponds to $D_3$ in Figures 2,5.
Another, somewhat related generalization of the induced matching lemma is in  \cite{GYtr} for so called transitive edge-colorings (where $aba$-free
corresponds to antichain coloring). It is worth observing that Theorem \ref{spattern} does not hold for the pattern $abab$, there are $1$-factorizations of
$K_{2^k}$  without an $abab$ pattern. This pattern is not a $2$-pattern and not transitive either. It is not clear what further patterns make
Theorem \ref{spattern} true.

For extended $s$-configurations we were not able to prove the Tur\'an-type Theorem \ref{tetel1}, but we could
prove a weaker result: a Ramsey-type theorem, Theorem \ref{extendedramsey} in the next subsection.

\subsection{Ramsey theory on Steiner triples}

Theorem \ref{tetel1} can be applied to prove the existence of Ramsey numbers in Steiner triple systems (the study of this
problem was initiated in \cite{GGY}).
Note that this is also strongly related to the Ramsey variant of Conjecture \ref{bes} which was studied
in \cite{ST}. However, in \cite{ST} the uniformity  is at least 4 unlike in our paper.

Here for historical reasons some of the terminology is a bit different from the one used for hypergraphs.
A {\em Steiner triple system of order $n$}, STS$(n)$, is an $n$-element set $V$,
called {\em points} and a set ${\cal{B}}$ of 3-element subsets of $V$ called
{\em blocks}, such that each pair of elements of $V$ appear in exactly one block
of ${\cal{B}}$. It is well-known that STS$(n)$ exists if and only if $n\equiv 1$ or $n\equiv 3$ $\pmod 6$, such values of $n$ are called {\em admissible}.

If {\em at most one} block covers each pair of elements of $V$, then we get  {\em partial} Steiner triple
systems, PTS$(n)$, this is the same as a linear 3-uniform
hypergraph and we call them {\em configurations}. As is customary, we assume that every point of a PTS$(n)$ is in at least
one block. The number of blocks containing $v\in V$ is called the {\em degree} of $v$.
A configuration $C$ is {\em unavoidable} if there is an $n_0=n_0(C)$ such that
every STS$(n)$ with $n\geq n_0$ must contain $C$. It is known \cite{BR,CR} that
all but two configurations with at most 4 blocks are unavoidable. The two
exceptions are $C_{14}$ and $C_{16}$ in Figure 6, the latter is called the {\em Pasch configuration}.\\

\begin{figure}[ht]
\centering
\includegraphics[scale=.7]{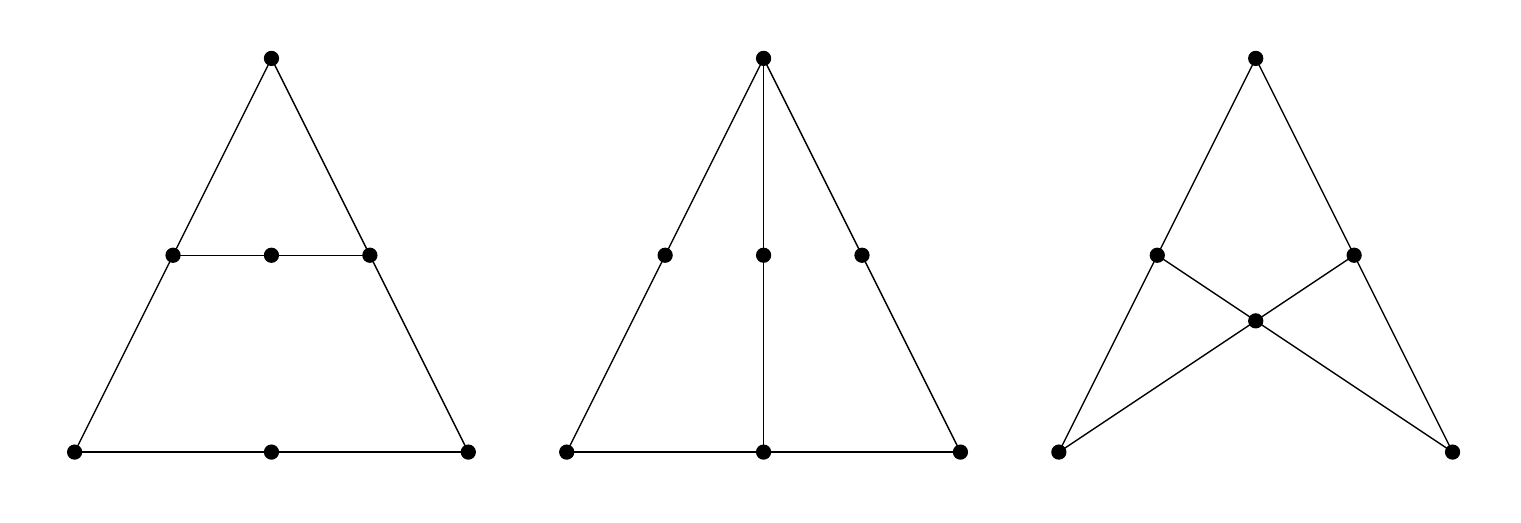}
\caption{Configurations $C_{14}$, $C_{15}$ (sail) and $C_{16}$ (Pasch configuration)}
\end{figure}

In \cite{GGY} a configuration was called {\em $t$-Ramsey} if for all large enough admissible
$n$ ($n\geq n_0(C, t)$), in every $t$-coloring of the blocks of any STS$(n)$
there is a monochromatic copy of $C$. If $C$ is $t$-Ramsey then
the smallest possible value of $n_0(C,t)$ is denoted by $R(C,t)$. Clearly,
a configuration is 1-Ramsey if and only if it is unavoidable.

A configuration $C=(V, {\cal{B}}) $ is called {\em acyclic} if
either $|{\cal{B}}|=1$, or it can be obtained from an acyclic
configuration $C'$ by adding a new block that intersects $V(C')$ in
at most one point. An equivalent definition is that $C=(V,
{\cal{B}}) $ does not contain {\em $i$-cycle} for any $i\ge 3$:
blocks $\{0,1,2\}+2j \pmod{2i}$ for $j=0,1,\dots,i-1$.

A configuration $C$ is called {\em graph-like} if every block
contains a point of degree one. In \cite{GGY} the following result
was proved.

\begin{theorem}[\cite{GGY}]\label{graph}
Acyclic configurations and graph-like configurations are both $t$-Ramsey
for any $t\geq 1$.
\end{theorem}

Note that an $s$-configuration is also a configuration, i.e. a partial triple system. In a $t$-coloring of the blocks of any STS$(n)$ the most
frequent color contains at least ${n-1\over 6}{n\over t}$ blocks. Thus Theorem \ref{tetel1} has the following corollary.
\begin{corollary}\label{cor}
Any $s$-configuration $C$ is $t$-Ramsey for all $s, t\geq 1$.
\end{corollary}

It is worth mentioning that Corollary \ref{cor} extends the ``acyclic'' part of Theorem \ref{graph} because of the following proposition.

\begin{proposition}\label{acext} Acyclic configurations are $s$-configurations for some $s$.
\end{proposition}

Notice that Corollary \ref{cor} does not follow for extended $s$-configurations. Nevertheless, we can prove the Ramsey property for them as well.

\begin{theorem}\label{extendedramsey} Any extended $s$-configuration $C$ is $t$-Ramsey for all $s, t\geq 1$.
\end{theorem}

In \cite{GGY} Ramsey properties of configurations with at most 4 blocks were studied.
The {\em sail} $C_{15}$ (sometimes also called a {\em fan}) is a configuration with blocks 123, 345, 561 and 147, see in Figure 6.  In \cite{GGY}
the following result was proved.
\begin{theorem}[\cite{GGY}]\label{andras}
All unavoidable configurations with at most 4 blocks, except possibly the sail, are $t$-Ramsey
for any $t\geq 1$.
\end{theorem}

Corollary \ref{cor} and Theorem \ref{extendedramsey} with some extra work allow us to extend Theorem \ref{andras} for configurations with at most 5 blocks.

\begin{theorem}\label{tetel2}
All unavoidable configurations with at most 5 blocks, except possibly the ones containing the sail, are $t$-Ramsey
for any $t\geq 1$.
\end{theorem}

The most interesting new $t$-Ramsey configuration in Theorem
\ref{tetel2} is the {\em wicket}, $D_4$ (see Figures 2 and 4), with blocks formed by the three rows and two columns of a
$3\times 3$ matrix of 9 points. In particular, the wicket is
$1$-Ramsey, i.e. unavoidable. But it has a much stronger property:
\begin{proposition}\label{strongunav} Every Steiner triple system except the Fano plane contains a wicket.
\end{proposition}

Note that it would be difficult to move to configurations with six
blocks since among them there is the {\em grid} (defined in
Subsection \ref{subturan} and conjectured to be avoidable
\cite{FR}).  But up to five blocks the crucial problem is the status
of the sail (see \cite{GGY}). In fact, the sail cannot be forced by density: for
$n\equiv 0$ $\pmod 3$
 there are sail-free PTS$(n)$'s with ${n^2\over 9}$ blocks. This is best possible as proved in \cite{FGY}.

In the next section we provide the tools.
Then in Section 3 we prove Theorems \ref{tetel1}, \ref{spattern}  and in Section 4 we prove Theorems \ref{extendedramsey} and \ref{tetel2}. In Section \ref{propproofs} we prove the two (easy but important) Propositions,  \ref{acext}, \ref{strongunav}.

\section{Tools}

For basic graph concepts see the monograph of Bollob\'as \cite{B}.\\
$V(G)$ and $E(G)$ denote the vertex-set and the edge-set
of the graph $G$. For $A\cap B=\emptyset$, $(A,B,E)$ denotes a bipartite graph $G=(V,E)$,
where $V=A\cup B$, and $E\subset A\times B$.
For a graph $G$ and a subset $U$ of its vertices,
$G|_U$ is the restriction of $G$ to $U$ .
$N(v)$ is the set of neighbors of $v\in V$.
Hence $|N(v)|=deg(v)=deg_G(v)$, the degree of $v$.
When $A,B$ are subsets of $V(G)$,
we denote by $e(A,B)$ the number of edges of $G$
with one endpoint in $A$ and the other in $B$.
For non-empty $A$ and $B$,$$d(A,B)=\frac{e(A,B)}{|A||B|}$$
is the {\em density} of the graph between $A$ and $B$.

\begin{definition}
The bipartite graph $G=(A,B,E)$ is $\eps$-regular if
$$X\subset A,\ Y\subset B,\ |X|>\eps|A|,\ |Y|>\eps|B|
\; \; imply \; \; |d(X,Y)-d(A,B)|<\eps,$$
otherwise it is $\eps$-irregular.
\end{definition}
We will often say simply that ``the pair $(A,B)$ is $\eps$-regular''
with the graph $G$ implicit.

In the proof of Theorem \ref{tetel1} the Regularity Lemma \cite{Sz}
plays a central role.
Here we will use the following variation of the lemma (see \cite{KS}).

\begin{lemma}[Regularity Lemma -- Degree form]\label{dreg}
For every $\eps>0$ and every integer $m_0$ there is an
$M_0=M_0(\eps, m_0)$ such that if $G=(V,E)$ is any graph on at
least $M_0$ vertices and $d \in[0,1]$ is any real number,
then there is a partition of the vertex-set $V$ into $l+1$ sets
(so-called clusters) $V_0,V_1,...,V_l$, and there is a subgraph
$G'=(V,E')$
with the following properties:
\begin{itemize}
\item $m_0\leq l\leq M_0$, \item $|V_0|\leq\eps|V|$, \item all
clusters $V_i,\,i\geq1,$ are of the same size $L$, \item
$deg_{G'}(v)>deg_G(v)-(d +\eps)|V|\; \text{for all} \; v\in V$,
\item $G'|_{V_i}=\emptyset$ ($V_i$ are independent in $G'$), \item
all pairs $G'|_{V_i\times V_j},\,1\leq i<j\leq l$, are
$\eps$-regular, each with a density 0 or exceeding $d$.
\end{itemize}
\end{lemma}

This form can be easily obtained by applying
the original Regularity Lemma (with a smaller value of $\eps$),
adding to the exceptional set $V_0$
all clusters incident to many irregular pairs,
and then deleting all edges between any other clusters
where the edges either do not form a regular pair
or they do but with a density at most $d$.

We will also use the following lemma.

\begin{lemma}\label{matching}
For every $\delta,c > 0$ there are positive constants $\eta, n_0$
with the following properties. Let $G$ be a graph on $n\geq n_0$ vertices
with at least $\delta n^2$ edges that is the edge disjoint union of matchings
$M_1, M_2, \ldots, M_m$ where $m\leq c n$. Then there exist an
$1\leq i \leq m$ and $A, B \subset V(M_i)$ (so the sets $A$ and $B$ are covered by $M_i$)
such that
\begin{itemize}
\item $|A|=|B|\geq \eta n$,
\item
$\left| E\left( G|_{A\times B}\right)\right| \geq \frac{\delta}{4} |A||B|$,
\item There are no $M_i$ edges within $A$, within $B$ and between $A$ and $B$.
\end{itemize}
\end{lemma}

\noindent
{\bf Proof of Lemma \ref{matching}:}
This lemma is basically identical to Lemma 2 in \cite{SS}.
It is also very similar to an argument in \cite{KS} (see the proof
of Theorem 3.2 in \cite{KS}). For the sake of completeness
we give the proof here.

Let $\eps > 0$ be small enough compared to $\delta$ and $c$.
We apply the degree form of the Regularity Lemma (Lemma \ref{dreg}) with parameters $\eps$ and
$d=\delta/2$, and let $G''$ denote the graph we get after removing
$V_0$ from $G'$. We still have
\beq\label{G''}|E(G'')|\geq \frac{\delta}{4} n^2.\eeq
Then using (\ref{G''}) there exists a matching $M_i$ for some
$1\leq i \leq m$, such that for $M = M_i \cap E(G'')$ we have
\beq\label{M}|M| \geq \frac{\delta}{4c} n.\eeq
Put $U=V(M)$ for the vertex set of $M$ and $U_i=V_i\cap U$.
Define
$$I = \{ i \; | \; |U_i| > 2 \eps |V_i| \}$$ and set
$L=\cup_{i\in I} U_i$ and $S=U\setminus L$.
Clearly, $|S|\leq 2\eps n$. Hence from (\ref{M}) we get
$|L|> \frac{|U|}{2}$ (since $\eps$ is small compared to $\delta$ and $c$),
and thus there exist two vertices $u, v \in L$
adjacent in $M$. Let $u\in V_i$ and $v\in V_j$.
Since there is an edge between $V_i$ and $V_j$ in $G''$, we must have a density of at least
$d=\frac{\delta}{2}$ between them. We may clearly select $A\subset U_i$
and $B\subset U_j$ so that
\beq\label{A,B}|A|\geq \eps |V_i| \; \mbox{and} \; |B|\geq \eps |V_j|\eeq
and that there is no edge of $M$ between $A$ and $B$.
By construction, there is no $M$ edge within $A$ and within $B$ since we
removed all edges within the clusters to get $G''$ and $A$ (and $B$) comes from one cluster.
Furthermore, (\ref{A,B}) and $\eps$-regularity imply that there is a density
of at least $\frac{\delta}{2} - \eps \geq \frac{\delta}{4}$ between $A$ and $B$,
as required. \qed

\section{Proof of Theorems \ref{tetel1} and \ref{spattern} }

{\bf Proof of Theorem \ref{tetel1} from Theorem \ref{spattern}:}
Using the well-known result of Erd\H{o}s and
Kleitman \cite{EK} (see also on page 1300 in \cite{GGL}) we find
a $3$-partite sub-hypergraph ${\cal{G'}}$ of ${\cal{G}}$ with at least
$$\frac{3!\delta}{3^3} n^2$$
edges.
Let $X_1, X_2, X_3$ be the vertex classes of this $3$-partite
hypergraph ${\cal{G'}}$.
For each $v\in X_1$ define the matching $M_v$ defined by $v$ between
$X_2$ and $X_3$ such that $(u,w)$ is an edge of $M_v$ if and only if
$(u,v,w)$ is a triple in ${\cal{G'}}$. $M_v$ is indeed a matching,
since ${\cal{G}}$ is a linear hypergraph.

Let $G$ be the properly colored bipartite graph that is the union
of the matchings $M_v$, where each $M_v$ is monochromatic.
Let $H$ be an $s$-pattern of ${\cal{H}}$.
Applying Theorem \ref{spattern}, $G$ contains $H$,
where a matching in $H$ goes to a matching in $G$ as well.
Adding back the corresponding vertices $v$,
we get ${\cal{H}}$ as a sub-hypergraph. \qed

Theorem \ref{spattern} in turn follows immediately from the following stronger lemma.
We prove this stronger statement because that is what we need in the proof of Theorem \ref{extendedramsey}.
If $H$ is an $s$-pattern,
then it is a bipartite graph (since it is acyclic), so we may assume that it has
a bipartition $V(H)=V_1\cup V_2$ with all edges going between $V_1$ and $V_2$.

\begin{lemma}\label{s-free}
For every $\delta,c >0$, integer $s\geq 1$ and $s$-pattern $H$ with bipartition $V(H)=V_1\cup V_2$,
there are positive constants $\gamma, n_0$ with the following property. Let $G$ be a
bipartite graph on $n\geq n_0$ vertices with bipartition $V(G)=U_1\cup U_2$ and with at least $\delta n^2$ edges
between $U_1$ and $U_2$ that is properly colored by at most $cn$ colors. Then $G$
contains at least $\gamma n$ vertex disjoint copies of $H$, where in the different copies of $H$ the same
matching always gets the same color in $G$ and $V_i$ is always embedded into $U_i$, $i=1, 2$.
\end{lemma}

\noindent
{\bf Proof of Lemma \ref{s-free}:}
Let the graph $G$ be the union
of matchings $$M_1, \ldots , M_m, m\leq cn,$$ where
each matching is the set of edges colored with a given color.
Assume that the number of vertices $n$ is sufficiently large.
Denote the
matchings in the definition of the $s$-pattern $H$ by
$M^1, \ldots , M^s$. Put $H_i = \cup_{j\geq i} M^j$.
We will define a sequence of subgraphs $G_1=G, G_2, \ldots , G_s$ of $G$.
Then we will show inductively from $s$ down to 1 that $G_i$ contains
at least $\gamma n$ vertex disjoint copies of $H_i$ with the properties claimed in
the lemma if $\gamma$ is sufficiently small. For $i=1$ we get the desired statement.

We define this sequence of subgraphs in the following way.
Put $G_1=G$. Let us apply Lemma \ref{matching}
for $G_1=G$ with
$$\delta_1 = \delta \;\; \mbox{and}\; \;  c_1=c.$$
We get $M_1'$ and $A_1$, $B_1$ as provided by Lemma \ref{matching}. Let $G_2 =  G_1|_{A_1\times B_1}$. Then
\begin{itemize}
\item $|A_1|=|B_1|\geq \eta_1 n$,
\item
$\left| E\left( G_2\right)\right| \geq \frac{\delta_1}{4} |A_1||B_1|$,
\item There are no $M_1'$ edges within $A_1$, within $B_1$ and between $A_1$ and $B_1$.
\end{itemize}
Put
$$\delta_2 = \frac{\delta_1}{16} \;\; \mbox{and}\; \;  c_2=\frac{c_1}{2\eta_1}.$$

Then indeed,
$$|E(G_2)| \geq \frac{\delta_1}{4} |A_1||B_1| = \frac{\delta_1}{16} (2|A_1|)(2|B_1|) = \delta_2 |V(G_2)|^2,$$
and
$$c_1 n = \frac{c_1}{2\eta_1} 2 \eta_1 n \leq \frac{c_1}{2\eta_1} 2 |A_1| = c_2 |V(G_2)|.$$

We continue in this fashion.
If $G_i, \delta_i, c_i, 1\leq i <s$ are already defined, then we apply Lemma \ref{matching}
for $G_i$ with $\delta_i$ and $c_i$ and we get $M_i'$ and $A_i$, $B_i$ as provided by Lemma \ref{matching}.
Let $G_{i+1}=G_i|_{A_i\times B_i}$. Then
\begin{itemize}
\item $|A_i|=|B_i|\geq \eta_i (|A_{i-1}|+|B_{i-1}|)$,
\item
$\left| E\left( G_{i+1}\right)\right| \geq \frac{\delta_i}{4} |A_i||B_i|$,
\item There are no $M_i'$ edges within $A_i$, within $B_i$ and between $A_i$ and $B_i$.
\end{itemize}
Put
$$\delta_{i+1} = \frac{\delta_i}{16} \;\; \mbox{and}\; \;  c_{i+1}=\frac{c_i}{2\eta_i}.$$

We will prove inductively from $s$ down to 1 that $G_i$ contains at
least $\gamma n$ vertex disjoint copies of $H_i$ with the properties claimed in
the lemma. For the base case $i=s$, $H_s=M^s$ is just a matching.
By the construction in $G_s=G_{s-1}|_{A_{s-1}\times B_{s-1}}$ we have
at least
$$\frac{\delta_{s-1}}{4} |A_{s-1}||B_{s-1}| \geq \frac{\delta}{4(16)^{s-2}} \left((\prod_{i=1}^{s-1} \eta_i) n\right)^2$$
edges, so we can find a monochromatic matching of size at least
$\gamma n |V(H)|/2$ in $G_s$, if
$$\frac{\delta \left((\prod_{i=1}^{s-1} \eta_i) n\right)^2}{4(16)^{s-2}c n} \ge \gamma n {|V(H)|\over 2},$$
i.e.
$$\gamma\le \frac{\delta (\prod_{i=1}^{s-1} \eta_i)^2}{2(16)^{s-2}c |V(H)|}.$$
Then we have at least $\gamma n$ monochromatic vertex disjoint copies
of $H_s=M^s$ in $G_s$, as desired.

Assume inductively that for some $i, 1\leq i< s$, $G_{i+1}$ contains at
least $\gamma n$ vertex disjoint copies of $H_{i+1}$ with the properties claimed in
the lemma. Note that from the definition of an $s$-pattern, every edge in $M^i$ either
has no edges of $H_{i+1}$ in its connected component, or touches an edge of $H_{i+1}$.
Therefore, $H_i$ is obtained from $H_{i+1}$ by attaching edges of $M^i$ to some vertices of
$H_{i+1}$ plus adding some isolated edges of $M^i$. Thus to get the desired copies of $H_i$ in $G_i$,
we add back the missing edges of $M^i$ to each copy of $H_{i+1}$ from the matching $M_i'$ in $G$,
so they all get the same color in $G$. This can be
done because each vertex of $A_i$ and $B_i$ is incident to exactly one edge from $M_i'$ which goes out of $A_i\cup B_i$.
Thus these edges are available to define $M^i$ for each copy of $H_{i+1}$, provided that
$|A_i|=|B_i|$ is large enough to allow to add all edges of $M^i$ in all of the at least $\gamma n$ copies.
This can be ensured by $$\gamma n {|V(H)|\over 2}\le (\prod_{j=1}^i \eta_j) n\le |A_i|=|B_i|,$$
i.e. by $\gamma\le {2(\prod_{j=1}^i \eta_j) \over |V(H)|}$.

Finally we get the desired statement from $i=1$. \qed

\section{Proof of Theorems \ref{extendedramsey} and \ref{tetel2}}

\noindent
{\bf Proof of Theorem \ref{extendedramsey}.} Assume that $C$ is an extended $s$-configuration defined by the extended $s$-pattern $H$. Let the bipartition of $H$ be $V(H)=V_1\cup V_2$ and
denote the matchings in the definition of the extended $s$-pattern $H$ by
$M^1, \ldots , M^s,M^*$. Let $uv$ be the single edge $M^*$ with $u\in V_1$ and $v\in V_2$
and let $H'$ be the disconnected $s$-pattern
resulting from $H$ after removing $uv$. In $H'$ we have two disconnected subgraphs
$H_1'$ and $H_2'$ such that $u\in V(H_1')$ and $v\in V(H_2')$.

Assume that $n$ is a sufficiently large admissible integer and we have a $t$-coloring of the blocks of an STS$(n)=(V,{\cal{B}})$. Partition $V$ into two almost equal parts $U_1$ and $U_2$.
Consider the complete bipartite graph $G$ between $U_1$ and $U_2$.  The edges of $G$ can be naturally colored with $t$ colors by assigning to $xy\in E(G)$ the color of the unique block $xyz$. We refer to this coloring as the {\em primary} coloring. On the other hand, there is also a natural proper coloring of the edges of $G$ with at most $n$ colors by assigning to $xy\in E(G)$ the vertex $z=f(x,y)$ of the unique block $xyz$. This coloring is called the $M$-coloring.

Let $G_1$ denote the subgraph of $G$ defined by the edges of the most frequent primary color, say red. Then $|E(G_1)| \geq \frac{n^2-1}{4t}\geq \frac{n^2}{5t}$.
Applying Lemma \ref{s-free} for $H'$ and $G_1$ with
$$\delta_1 = \frac{1}{5t} \;\; \mbox{and}\; \;  c_1=1,$$
we find at least $\gamma_1 n$ vertex disjoint red copies of $H'$ in $G_1$, where in the different copies of $H'$ the same
matching always gets the same $M$-color in $G_1$ and $V_i$ is always embedded into $U_i$, $i=1, 2$.
Let $A_1\subset U_1$ denote the set of embedded images of $u$ in these copies of $H'$, similarly let $B_1\subset U_2$
denote the set of embedded images of $v$ in these copies of $H'$. Let $L_1,\dots,L_m$ denote the copies of $H'_1$ and $u_i, 1\leq i \leq m$ the image of $u$ in these copies. Similarly, $R_1,\dots,R_m$ denote the copies of $H'_2$ and $v_i, 1\leq i \leq m$ the image of $v$ in these copies.
So $A_1 = \{u_1,...,u_m\}$ and $B_1 = \{v_1,...,v_m\}$, and the $i$-th copy of $H'$ consists of $L_i$ and $R_i$.

A copy of $H'$, $L_i\cup R_i$,  is bad if for some edge $xy$ in any of the copies $L_j\cup R_j$ of $H'$, $f(x,y)$ is a vertex in $L_i\cup R_i$,  otherwise $L_i\cup R_i$ is good.
Note that there are at most $s$ bad copies since there are at most $s$ $M$-colors on the edges of $H'$. Let us remove the bad copies and denote by $A_1'$ and $B_1'$ the set of remaining vertices in $A_1$ and $B_1$.
Then
\beq\label{A_1'}|A_1'|=|B_1'|\geq |A_1|-s\geq \gamma_1 n - s \geq \frac{\gamma_1}{2} n.\eeq
Assume that $L_i,R_j$ are good copies and the primary color of $u_iv_j$ is red. In this case we try to extend the $s$-pattern with the edge $u_iv_j$ and find a red copy of $C$. We have to avoid the following two exceptional situations for success.

\begin{itemize}

\item $f(u_i,v_j)=f(a,b)$ for some edge $ab$ in $L_i$ or in $R_j$. In this case we cannot extend the $s$-pattern with the edge $u_iv_j$ because its $M$-color is not a new $M$-color. There are at most $s|A_1'|$ possibilities for this situation. Indeed, the number of $M$-colors in $H'$ is $s$ and in each of these $M$-colors there can be at most $|A_1'|$ edges between $A_1'$ and $B_1'$.
\item $f(u_i,v_j)$ is covered by $V(L_i)\cup V(R_j)$. In this case we can extend the $s$-pattern with the edge $u_iv_j$ to the required extended $s$-pattern but the point $f(u_i,v_j)$ of the red block  $u_iv_jf(u_i,v_j)$  may not be well placed, namely, it will not be a new vertex. However, each fixed $M$-color can appear at most once (since the copies $L_i,R_j$ are disjoint), thus at most $n$ pairs $u_iv_j$ can be in this situation. Indeed, for each $z$ in $L_i$, the pair $z,u_i$ is only in one block, hence there can be at most one $j$ such that $z = f(u_i,v_j)$, and similarly for $z$ in $R_j$. So each of the at most $n$ vertices in $L_1,...,L_m,R_1,...,R_m$ ``ruins" at most one pair $u_i,v_j$.
\end{itemize}

Then if we have a non-exceptional red $u_iv_j$, this will
be the image of $uv$. To get a red copy of $H$ in $G_1$, we take $L_i$, the copy
of $H_1'$ containing $u_i$ and $R_j$, the copy of $H_2'$ containing $v_j$.
Finally, adding back the corresponding 3rd vertices of the blocks (which are
disjoint from this red copy of $H$ by construction) we get a red copy of $C$.

Thus we may assume that there is no such non-exceptional red $u_iv_j$.
Then in $G|_{A_1'\times B_1'}$ apart from at most $s|A_1'|+n$ exceptional edges,
all edges are colored (in the primary coloring) with the remaining $(t-1)$ colors
(other than red). Let $G_2$
denote the subgraph of the most frequent color out of these $(t-1)$ colors in $G|_{A_1'\times B_1'}$. Then
$$|E(G_2)| \geq \frac{|A_1'||B_1'| - s |A_1'|-n}{t-1} \geq \frac{1}{2(t-1)}|A_1'||B_1'|.$$
Indeed, this follows from
$$s|A_1'|+ n \leq \frac{|A_1'||B_1'|}{2},$$
which in turn follows from
$$2s+\frac{4}{\gamma_1} \leq |B_1'|$$
(using (\ref{A_1'}) and the fact that $n$ is sufficiently large).

We will apply Lemma \ref{s-free} for $H'$ and $G_2$ with
$$\delta_2 = \frac{1}{8(t-1)} \;\; \mbox{and}\; \;  c_2=\frac{c_1}{\gamma_1}.$$
Then indeed,
$$|E(G_2)| \geq \frac{1}{2(t-1)} |A_1'||B_1'| = \frac{\delta_1}{8(t-1)} (2|A_1'|)(2|B_1'|) = \delta_2 |V(G_2)|^2,$$
and
$$c_1 n = \frac{c_1}{\gamma_1} \gamma_1 n \leq \frac{c_1}{\gamma_1} 2 |A_1'| = c_2 |V(G_2)|.$$
Applying Lemma \ref{s-free}
we find at least $\gamma_2 |V(G_2)|$ vertex disjoint copies of $H'$ in $G_2$, where in the different copies of $H'$ the same
matching always gets the same $M$-color in $G_2$ and $V_i$ is always embedded into $U_i$, $i=1, 2$.

We continue in this fashion. We will apply Lemma \ref{s-free} for $H'$ and $G_i$ with
$$\delta_i = \frac{1}{8(t-(i-1))} \;\; \mbox{and}\; \;  c_i=\frac{c_{i-1}}{\gamma_{i-1}}.$$
Then indeed,
$$|E(G_i)| \geq \delta_i |V(G_i)|^2,$$
and
$$c_{i-1} |V(G_{i-1})| = \frac{c_{i-1}}{\gamma_{i-1}} \gamma_{i-1} |V(G_{i-1})| \leq \frac{c_{i-1}}{\gamma_{i-1}} 2 |A_{i-1}'| = c_i |V(G_i)|.$$
Applying Lemma \ref{s-free}
we find at least $\gamma_i |V(G_i)|$ vertex disjoint copies of $H'$ in $G_i$, where in the different copies of $H'$ the same
matching always gets the same $M$-color in $G_i$ and $V_i$ is always embedded into $U_i$, $i=1, 2$.
Note that all edges of $G_i$ have the same primary color, so these copies of $H'$ we find are monochromatic in this color.
Let $A_i\subset A_{i-1}'$ denote the set of embedded images of $u$ in these copies of $H'$, similarly let $B_i\subset B_{i-1}'$
denote the set of embedded images of $v$ in these copies of $H'$.
We remove the bad copies of $H'$ and denote by $A_i'$ and $B_i'$ the set of remaining vertices in $A_i$ and $B_i$.

Again if there is a non-exceptional edge in $G_i$ between $A_i'$ and $B_i'$, then we are done.
Otherwise in $G|_{A_i'\times B_i'}$ apart from at most $is|A_i|+in$ exceptional edges,
all edges are colored (in the primary coloring) with the remaining $(t-i)$ colors. Let $G_{i+1}$
denote the subgraph of the most frequent color out of these $(t-i)$ colors in $G|_{A_i'\times B_i'}$. Then
$$|E(G_{i+1})| \geq \frac{|A_i'||B_i'| - i s |A_i'|-in}{t-i}\geq \frac{1}{2(t-i)}|A_i'||B_i'|,$$
(assuming that $n$ is sufficiently large).

Finally we arrive at $G_t$ between $A'_{t-1}$ and $B'_{t-1}$, where
all but $O(n)$ edges are of the last primary color, and in particular the density of $G_t$ is at least 1/2.
Applying Lemma \ref{s-free} for $H'$ and $G_t$ with
$$\delta_t = \frac{1}{2} \;\; \mbox{and}\; \;  c_t=\frac{c_{t-1}}{\gamma_{t-1}},$$
we find at least $\gamma_t |V(G_t)|$ vertex disjoint copies of $H'$ in $G_t$, where in the different copies of $H'$ the same
matching always gets the same $M$-color in $G_t$ and $V_i$ is always embedded into $U_i$, $i=1, 2$.
Let $A_t\subset A_{t-1}'$ denote the set of embedded images of $u$ in these copies of $H'$, similarly let $B_t\subset B_{t-1}'$
denote the set of embedded images of $v$ in these copies of $H'$.
We remove the bad copies of $H'$ and denote by $A_t'$ and $B_t'$ the set of remaining vertices in $A_t$ and $B_t$.
Now there must be a non-exceptional edge in $G_t$ between $A_t'$ and $B_t'$ because this is the last color available.
Indeed, other than at most $ts|A_t'|+tn$ exceptional edges,
all edges must have this color. Then we are done similarly as before. \qed

\medskip

\noindent
{\bf Proof of Theorem \ref{tetel2}.} It is known that there are 56 non-isomorphic configurations with 5 blocks (see \cite{CR}).
Because of Theorem \ref{graph} and the assumptions of Theorem \ref{tetel2} we can exclude the following configurations:
\begin{itemize}
\item Acyclic and graph-like configurations,
\item Configurations containing avoidable configurations with 4 blocks, the Pasch configuration $C_{16}$  and $C_{14}$,
\item Configurations containing the sail, $C_{15}$.
\end{itemize}

Let $C$ be a configuration with 5 blocks different from the ones
listed above. Then $C$ must contain an $i$-cycle for some $i\ge 3$.
We distinguish between two cases. If a configuration is a 2- or
3-configuration or an extended 2-configuration, then we know by
Corollary \ref{cor} and Theorem \ref{extendedramsey} that it is
$t$-Ramsey.

\noindent
{\bf Case 1:} $C$ contains a triangle $T$ ($3$-cycle) with blocks $123,345,561$.\\
Assume that the other two blocks are $B_1,B_2$.  Note that both $B_i$'s intersect $T$ in at most one point because otherwise we would get a copy of $C_{14},C_{16}$ or $C_{15}$.

One of the $B_i$'s, say $B_1$ must intersect $T$ in a point of
degree one, say in $\{2\}$, because $C$ is not graph-like. Set
$P=\{5\}$, $Q=\{2\}$ and let $R$ be a degree one point of $B_2$ (it
exists since $B_2$ intersects both $B_1$ and $T$ in at most one
point). Define the ``pattern" obtained by labeling with $p,q,r$ the
pairs in the blocks containing $P,Q,R$, respectively. For example
$T\cup B_1$ defines the pattern $pqp,q$. When $R,Q$ or $R,P$ are in
the same block (second and fifth possibilities in Subcase 1.2 below)
we will not use color $r$ but only $q$ and $p$.

\noindent
{\bf Subcase 1.1:} $B_2$ does not intersect $T\cup B_1$.  ($(11,5$)-configuration.)\\
Now $C$ is a $3$-configuration (so it is $t$-Ramsey by Corollary \ref{cor})
based on the $3$-pattern with components $pqp,q$ and $r$. See Figure 7
for the $3$-configuration and Figure 8 for the resulting $3$-pattern. From now on for
simplicity we just show the patterns.

\begin{figure}[ht]
\centering
\includegraphics[scale=.8]{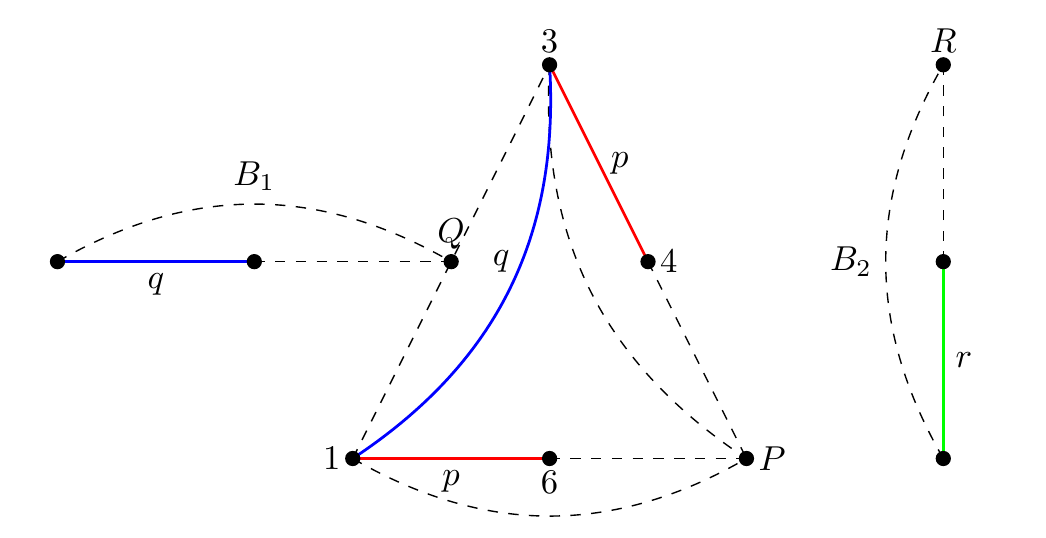}
\caption{The 3-configuration in Subcase 1.1}
\end{figure}

\begin{figure}[ht]
\centering
\includegraphics[scale=.9]{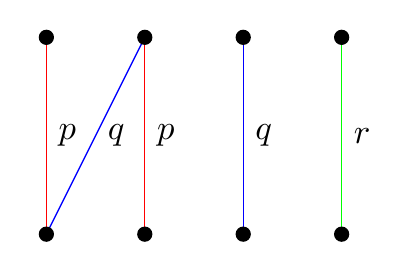}
\caption{The resulting 3-pattern in Subcase 1.1}
\end{figure}

\noindent
{\bf Subcase 1.2:} $B_2$  intersects $T\cup B_1$ in one point. ($(10,5)$-configurations.)\\
W.l.o.g. we may assume that the intersection point is one of $B_1
\setminus T$, $\{2\}=Q$, $\{3\}$, $\{4\}$, $\{5\}=P$ (or we get
isomorphic configurations). These configurations are all $2$- or
$3$-configurations. They are based on the 2- or 3-patterns $pqp,qr$;
$pqp,q,q$; the star with edges $p,q,r$ with a $p$ edge attached at
$q$ plus a $q$ component; $pqpr,q$ and $pqp,p,q$ (see Figure 9).

\begin{figure}[ht]
\centering
\includegraphics[scale=.9]{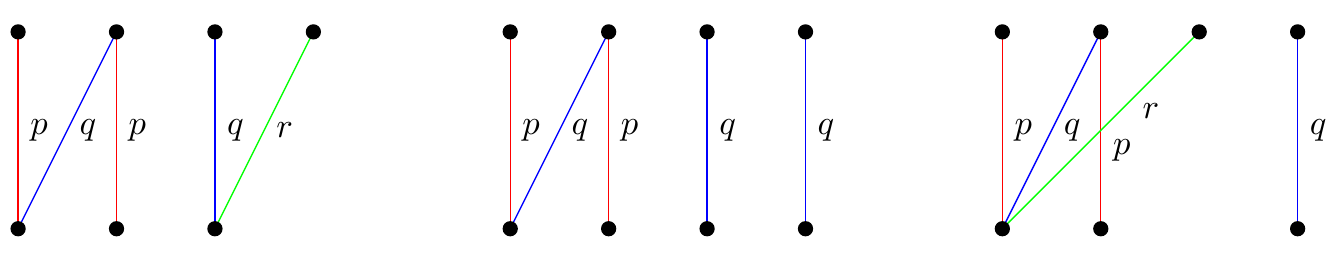}

\includegraphics[scale=.9]{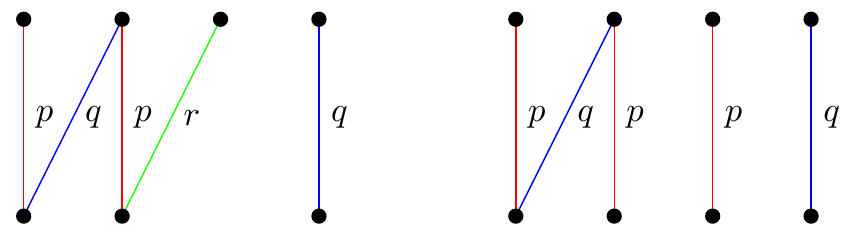}
\caption{The patterns in Subcase 1.2}
\end{figure}

\noindent
{\bf Subcase 1.3:} $B_2$  intersects $T\cup B_1$ in two points. ($(9,5)$-configurations.)\\
One intersection point must be in $B_1\setminus T$ and the other is w.l.o.g. at $\{3\}$, $\{4\}$ or $\{5\}$. The first is a $3$-configuration, $D_1$, based on the $3$-pattern $pqrq$ with a $p$ edge at the midpoint. The second is an extended 2-configuration, $D_3$, based on the extended 2-pattern $pqprq$
(in fact it is also a $3$-configuration, see Figure 5.) The third, $D_2$, is a $2$-configuration based on the $2$-pattern $pqp,pq$ (see Figure 10).

\begin{figure}[ht]
\centering
\includegraphics[scale=.9]{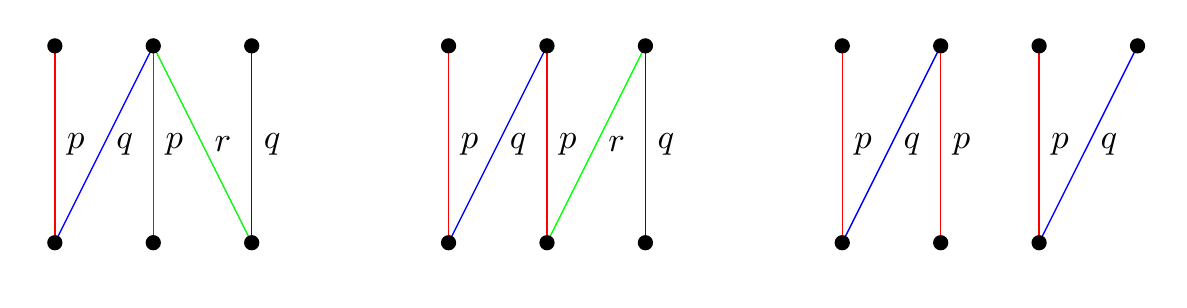}
\caption{The patterns in Subcase 1.3}
\end{figure}

\noindent
{\bf Case 2:} $C$ contains no triangle but contains a $4$-cycle with blocks $123$, $345$, $567$, $781$.\\
Then this 4-cycle can be extended (avoiding a graph-like configuration) in two ways.

\noindent
{\bf Subcase 2.1:} A new block intersects the $4$-cycle in one of its degree one points, say in $\{2\}$.  ($(10,5)$-configuration.)\\
This is a $3$-configuration based on the $3$-pattern $pqrp,r$.
See Figure 11
for the $3$-configuration and Figure 12 for the resulting $3$-pattern.

\begin{figure}[ht]
\centering
\includegraphics[scale=.8]{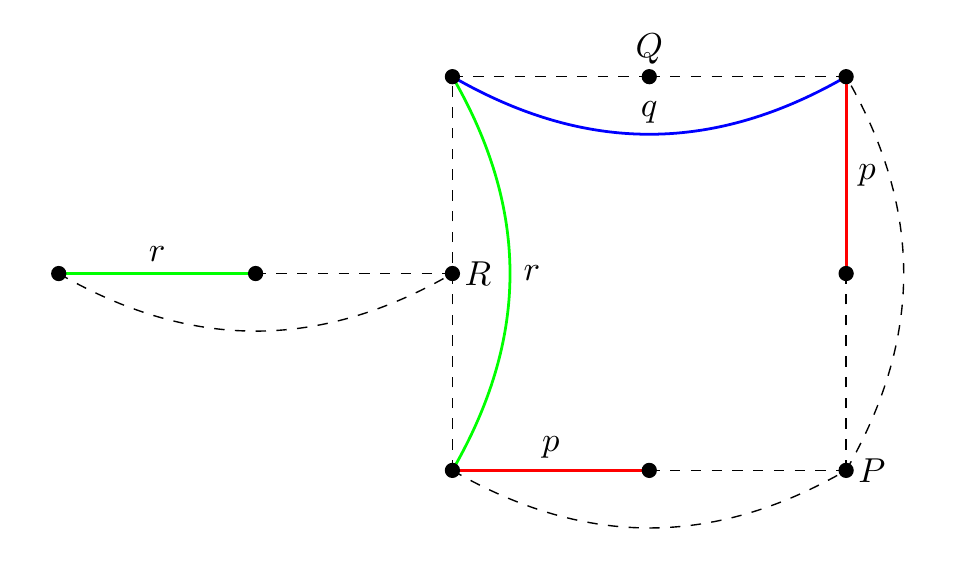}
\caption{The 3-configuration in Subcase 2.1}
\end{figure}

\begin{figure}[ht]
\centering
\includegraphics[scale=.9]{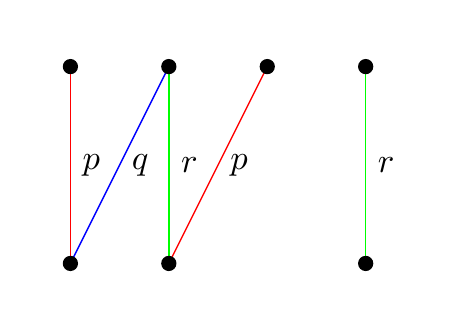}
\caption{The resulting 3-pattern in Subcase 2.1}
\end{figure}

\noindent
{\bf Subcase 2.2:}  A new block intersects the $4$-cycle in two opposite degree one points, say it is $269$. ($(9,5)$-configuration.)\\
 This is an extended $2$-configuration, $D_4$ (the wicket), based on the extended $2$-pattern $pqrpq$ (see Figure 3 where $a,b,c$ are used instead of $p,q,r$).
 Note that this is the only place where we really used the full power of Theorem \ref{extendedramsey}.

Cases 1,2 cover all possibilities because the $5$-cycle is graph-like. \qed


\section{Proof of Propositions \ref{acext} and \ref{strongunav}}\label{propproofs}

{\bf Proof of Proposition \ref{acext}.} Consider an acyclic
configuration $C$, it has a block $B=xyz$ such that $y,z$ are points
of degree one. By induction, the acyclic configuration $C'$ obtained
by the removal of $B$ is an $s$-configuration. If $x$ is an
augmenting point of the $s$-pattern $S$ of $C'$ then $x$ is also an
augmenting point of the $s$-pattern $S\cup \{y,z\}$. Otherwise $C$
is an $(s+1)$-configuration obtained from the $(s+1)$-pattern $S\cup
\{x,y\}$ with the new augmenting point $z$. \qed

\medskip
\noindent
{\bf Proof of Proposition \ref{strongunav}. } It is easy to see
(since STS$(7)$, the Fano plane, is unique) that for $n>7$ any
STS$(n)$ contains two disjoint blocks,
$X=\{x_1,x_2,x_3\},Y=\{y_1,y_2,y_3\}$. Consider the proper coloring
of the complete bipartite graph $[X,Y]$ where the color of $x_iy_j$
is defined by the third point of the block containing the pair
$(x_i,y_j)$. We claim that there is a ``rainbow'' matching in this
coloring, three disjoint pairs with three different colors. Indeed,
assume that $(x_1,y_1)$ has color $a$, then one of
$(x_2,y_2),(x_2,y_3)$, say $(x_2,y_2)$ has color $b\ne a$. Assume
w.l.o.g. that $(x_3,y_3)$ has color $a$ (otherwise we have the
rainbow matching). The pairs $(x_1,y_2),(x_2,y_1)$ must be colored
with the same color, different from both $a,b$, say color $c$. The
same is true for the pairs $(x_2,y_3),(x_3,y_2)$, they are colored
with $d$ which is different from $a,b,c$. But then
$(x_1,y_3),(x_2,y_1),(x_3,y_2)$ is a rainbow matching since
$(x_1,y_3)$ cannot have color $c$ or $d$ because it touches edges
with these colors. Then the blocks on these pairs together with
blocks $X,Y$ define a wicket. \qed

\bigskip

\noindent
{\bf Acknowledgement. } We appreciate the detailed careful comments of the referees, they have much improved the presentation.

\end{document}